\documentclass[conference]{IEEEtran}
\IEEEoverridecommandlockouts

\usepackage{textcomp}
\usepackage{dsfont,amssymb,amsmath,fancyhdr,mdframed}
\usepackage[dvipsnames]{xcolor}
\usepackage{comment}
\usepackage{dsfont}
\usepackage{graphicx}%
\usepackage{multicol,lipsum}
\usepackage{tikz}
\usepackage{hyperref}
\usepackage{subcaption}
\usepackage[prependcaption,colorinlistoftodos]{todonotes}
\usepackage{algorithm, algorithmic}
\usetikzlibrary{arrows.meta} 
\usetikzlibrary{fit, positioning}
\usetikzlibrary{calc}
\usetikzlibrary{decorations.pathmorphing} 

\tikzstyle{block} = [rectangle, minimum width=1cm, minimum height=1cm, text centered, draw=black]
\tikzstyle{tallblock} = [rectangle, minimum width=.5cm, minimum height=1cm, text centered, draw=black]
\tikzstyle{line} = [thick,-,>=stealth]
\tikzstyle{arrow} = [thick,->,>=stealth]
\tikzstyle{roundedblock} = [rectangle, minimum width=4cm, minimum height=2cm, text centered, draw=black, rounded corners=0.2cm]

\newcommand{\mc}{\mathcal}

\newcommand{\col}{\mathrm{col}}

\newcommand{\blkdiag}{\mathrm{blkdg}}
\newcommand{\blkmat}{\mathrm{blkmat}}

\newtheorem{theorem}{Theorem}

\newtheorem{lemma}{Lemma}

\def\BibTeX{{\rm B\kern-.05em{\sc i\kern-.025em b}\kern-.08em
    T\kern-.1667em\lower.7ex\hbox{E}\kern-.125emX}}
\begin{document}

\title{Fast Newton methods for linear-quadratic dynamic games with application to autonomous vehicle platooning and intersection crossing
\thanks{This work was partially supported by the Transport \& Mobility Institute (TMI), TU Delft, and by the ERC under project ARGON.}
}

\author{\IEEEauthorblockN{Reza Rahimi Baghbadorani}
\IEEEauthorblockA{\textit{Erasmus University Rotterdam} \\
\texttt{rahimibaghbadorani@rsm.nl}}
\and
\IEEEauthorblockN{Sergio Grammatico}
\IEEEauthorblockA{\textit{Delft University of Technology} \\ 
\texttt{s.grammatico@tudelft.nl}}
}

\maketitle

\begin{abstract}
We consider constrained linear–quadratic dynamic games arising in autonomous vehicle platooning, intersection crossing and other cooperative driving scenarios. Infinite-horizon Nash equilibria are reformulated as receding-horizon affine variational inequalities with special structure. Exploiting this formulation, we design Newton-type algorithms with local quadratic convergence. The resulting methods achieve extremely fast convergence, making them well suited for real-time and embedded receding-horizon control in safety-critical traffic applications. Simulations of platooning and intersection crossing demonstrate substantial performance gains over first-order and operator-splitting approaches, hence high application potential.

\end{abstract}


\section{Introduction}\label{sec: intro}

The rapid advancement of Intelligent Transportation Systems (ITS) requires sophisticated methodologies to manage the complex interactions in modern infrastructures, such as in traffic flow control and multi-vehicle routing \cite{talebpour2016influence,di2014distributed,spica2020real}. To address these challenges, dynamic game theory and automatic control offer robust frameworks for analyzing and designing the interconnections in multi-agent systems \cite{wang2021game,bacsar1998dynamic}.

\textit{Dynamic game theory} provides a formal framework for modeling feedback control where agents' objectives and constraints are coupled through shared system dynamics \cite{bacsar1998dynamic,basar1976uniqueness}. In these settings, each agent seeks to minimize an individual cost function while anticipating the strategic maneuvers of others. These interactions are frequently characterized by the {Open-Loop Nash Equilibrium (OL-NE)}, where control inputs are determined based on the initial state and the anticipated actions of other participants \cite{monti2024feedback}.

Remarkably, the OL-NE problem for linear-quadratic games can be reformulated as a (strongly) monotone affine Variational Inequality (VI) \cite{monti2024feedback,benenati2025linear}. Solving VIs efficiently is critical, especially for receding-horizon real-time control implementations, as they require computing a new equilibrium at every sampling interval~\cite{benenati2025linear}. Therefore, the computational method must be both fast and reliable. However, the complexity of these VIs increases rapidly with the number of agents, the length of the control horizon, and the density of constraints.

While recent literature has demonstrated the efficacy of first-order splitting methods, such as the Douglas-Rachford (DR) algorithm \cite{baghbadorani2025douglas} for solving affine VIs in OL-NE, there remains a gap in the application of higher-order information. This paper investigates whether Newton-type methods can provide the necessary acceleration for real-time applications. 


\subsubsection*{Contributions}

In our preliminary numerical experiments, Newton methods significantly outperform all others for strongly monotone VIs \cite{mignoni2025monviso}, e.g. forward-backward descent (FB), extragradient descent, Nesterov’s accelerated gradient descent (NAGD), projected reflected gradient descent, adaptive Golden ratio (aGRAAL), and Douglas--Rachford splitting method (DR). This motivates us to tailor Netwon methods for real-time control of autonomous systems.

\begin{itemize}
    \item {Algorithm development:} We propose a \textit{smoothed Netwon method} tailored to the structure of infinite-horizon OL-NE linear-quadratic games \cite{benenati2025linear}.
    \item {Theoretical guarantees:} We provide local superlinear convergence (Theorem \ref{Th: local convergence}) and global convergence via a line-search mechanism (Theorem \ref{thm: global conv}).
    \item {Numerical validation:} We demonstrate our approach in two ITS applications: vehicle platooning and autonomous intersection crossing, benchmarking its performance in receding-horizon control scenarios.
\end{itemize}

\subsubsection*{Notation}
We adopt a standard notation from \cite{benenati2025linear, baghbadorani2025douglas}.
For a closed convex set $\mathcal C$, the metric projection onto $\mathcal C$ is given by $\pi_{\mathcal C}(u) := \arg\min_{y \in \mathcal C} \|u - y\|$. An operator $F : \mathcal C \to \mathbb{R}^n$ is $L$-Lipschitz continuous if there exists $L > 0$ such that $\|F(u) - F(y)\| \le L \|u - y\|, \quad \forall u,y \in \mathcal C$. Moreover, $F$ is (strongly) monotone if there exists $\mu \ge 0$ (respectively $\mu > 0$) such that $\langle F(u) - F(y), u - y \rangle \ge \mu \|u - y\|^2, \quad \forall u,y \in \mathcal C$. For an affine operator $F(x) = Mx + q$ and a feasible set $\mathcal C$, we denote the associated variational inequality problem by $\mathrm{AVI}(\mathcal C, M, q)$. We denote $r_k$ as the residual in the iterative method for solving the VI, a metric to measure the optimality of the iterative method at the point $u_k$~\cite[Prop.~1.5.8]{facchinei2003finite}, defined as $r_k := \| u_k - \pi_\mathcal{C}(u_k - F(u_k)) \|$.
We denote by $\mathcal S_T^n$ a sequence of vectors in $\mathbb{R}^n$ of length $T \in \mathbb{N} \cup \{\infty\}$. For a sequence $w \in \mathcal S_T^n$, its element at time index $t \in \{0,\dots,T-1\}$ is written as $w[t]$. 
\section{Open-Loop Nash equilibrium via receding horizon Affine Variational Inequality} \label{sec: OL-NE as AVI]}


We consider a discrete-time linear system:
\begin{equation}
    x[t+1] = A x[t] + \sum_{i\in\mathcal I} B_i u_i[t], \quad x[0] = x_0,
\end{equation}
where \(x[t]\in\mathbb{R}^n\), \(u_i[t]\in\mathbb{R}^m\), and \(\mathcal I = \{1,\dots,N\}\). For a joint input sequence \(\mathbf u = (u_1, \dots, u_N)\), the resulting trajectory is denoted by \(\zeta(t,x_0,\mathbf u)\). Each agent \(i\) aims to minimize the infinite-horizon quadratic cost
\begin{equation}\label{eq: inf-LQgame}
    J_i^\infty(u_i, \mathbf u_{-i}, x_0) = \sum_{t=0}^{\infty} \frac12 x[t]^\top Q_i x[t] + \frac12 u_i[t]^\top R_i u_i[t],
\end{equation}
where \(Q_i \succeq 0\) and \(R_i \succ 0\). The agents are subject to affine state and input constraints, defined as
\begin{align*}
    \mathcal U_i(\mathbf u_{-i}[t]) &:= \Bigl\{ u_i[t] : \sum_{j\in\mathcal I} D_j u_j[t] + d_u \le 0 \Bigr\},\\
    \mathcal X &:= \{ x[t] : D_x x[t] + d_x \le 0 \},
\end{align*}
and for a finite horizon \(T\in\mathbb N\), we define the collective feasible set as
\begin{align*}
    \mathcal U_T(x_0) &:= \\
    &\Bigl\{ \mathbf u : u_i[t] \in \mathcal U_i(\mathbf u_{-i}[t]),\;\zeta(t,x_0,\mathbf u) \in \mathcal X,\; \forall i, t<T \Bigr\},
\end{align*}
which can equivalently be expressed as a set of affine inequalities $\mathcal U_T(x_0) = \{ \mathbf u : D \mathbf u + d \le 0 \}$.

An \emph{infinite-horizon open-loop Nash equilibrium (OL-NE)} is a feasible sequence \(\mathbf u^\star \in \mathcal U_\infty(x_0)\) such that the state converges to zero, $\lim_{t \to \infty} \zeta(t,x_0,\mathbf u^\star) = 0$, and no agent can reduce its own cost by a unilateral deviation, i.e., for all \(i \in \mathcal I\) and all \(u_i\) such that \((u_i, \mathbf u^\star_{-i}) \in \mathcal U_\infty(x_0)\),
\[
J_i^\infty(u_i^\star, \mathbf u^\star_{-i}, x_0) \le J_i^\infty(u_i, \mathbf u^\star_{-i}, x_0).
\]
Since this problem is generally intractable, we rely on a finite-horizon approximation based on stabilizing feedback equilibria \cite{benenati2025linear,monti2024feedback}. To approximate the infinite-horizon problem, we introduce the finite-horizon cost with terminal term
\begin{align*}
    &J_i^T(u_i, \mathbf u_{-i}, x_0) := \sum_{t=0}^{T-1} \frac12 \zeta(t,x_0,u_i,\mathbf u_{-i})^\top Q_i \zeta(t,x_0,u_i,\mathbf u_{-i}) \\
    &+ \frac12 u_i[t]^\top R_i u_i[t] + \frac12 \zeta(T,x_0,u_i,\mathbf u_{-i})^\top \hat P_i \zeta(T,x_0,u_i,\mathbf u_{-i}),
\end{align*}
where $\hat P_i$ is the terminal cost obtained from the augmented Riccati equation.  Specifically, the coupled Riccati equations
\begin{align*}
    P_i &= Q_i + A^\top P_i \Bigl(A + \sum_{j\in\mathcal I} B_j K_j\Bigr), \\
    K_i &= -R_i^{-1} B_i^\top P_i \Bigl(A + \sum_{j\in\mathcal I} B_j K_j\Bigr)
\end{align*}
define stabilizing feedback gains $K_i$ and closed-loop dynamics $A_{\mathrm{cl}} = A + \sum_j B_j K_j$. For each agent $i$, the augmented system for unilateral deviations is
\begin{align*}
    \hat x_i[t+1] &= \hat A_i \hat x_i[t] + \hat B_i u_i[t], \\
    \hat A_i &= \begin{bmatrix} A & \sum_{j\neq i} B_j K_j \\ 0 & A_{\mathrm{cl}} \end{bmatrix}, \quad
    \hat B_i = \begin{bmatrix} B_i \\ 0 \end{bmatrix},
\end{align*}
with stabilizing solution \(\hat P_i \succeq 0\) solving
\begin{align*}
    \hat P_i &= \hat Q_i + \hat A_i^\top \hat P_i (\hat A_i + \hat B_i \hat K_i), \,\, \hat Q_i = \mathrm{blkdiag}(Q_i,0)\\
    \hat K_i &= -R_i^{-1} \hat B_i^\top \hat P_i (\hat A_i + \hat B_i \hat K_i).
\end{align*}

Next, we reformulate the approximated infinite-horizon game \eqref{eq: inf-LQgame} as a variational inequality. By stacking the states and control inputs over the horizon $T$, we have
\begin{equation}
    \mathbf x = \begin{bmatrix} x[1] \\ \vdots \\ x[T] \end{bmatrix}, \qquad
    \mathbf u_i = \begin{bmatrix} u_i[0] \\ \vdots \\ u_i[T-1] \end{bmatrix},
\end{equation}
and therefore $\mathbf x = \Theta x_0 + \sum_{i\in\mathcal I} \Gamma_i \mathbf u_i$, where
\begin{equation}
    \Theta = \begin{bmatrix} A \\ A^2 \\ \vdots \\ A^T \end{bmatrix}, \quad
    \Gamma_i =
    \begin{bmatrix}
    B_i & 0 & \cdots & 0 \\
    AB_i & B_i & \cdots & 0 \\
    \vdots & \vdots & \ddots & \vdots \\
    A^{T-1}B_i & A^{T-2}B_i & \cdots & B_i
    \end{bmatrix}.
\end{equation}
Now, Let $\bar R_i = I_T \otimes R_i$ and $\bar Q_i = \mathrm{blkdiag}(I_{T-1} \otimes Q_i, \hat P_i)$. Then the stacked gradient of \(J_i^T\) is
\begin{equation}
    \nabla_{\mathbf u_i} J_i^T(\mathbf u) = \bar R_i \mathbf u_i + \Gamma_i^\top \bar Q_i \Bigl( \Theta x_0 + \sum_{j\in\mathcal I} \Gamma_j \mathbf u_j \Bigr).
\end{equation}
Next, we define operator $F(\mathbf u) := M \mathbf u + q$, where
\begin{align*}
    M &:= \blkdiag(\bar{R}_i)_{i\in\mc I} + \blkmat(\Gamma_i^{\top}\bar{Q}_i\Gamma_j)_{(i,j)\in\mc I^2},\\ 
    q &:= \col(\Gamma_i^{\top}\bar{Q}_i\Theta x_0)_{i\in\mc I}.
\end{align*}
The following result establishes a connection between infinite-horizon and finite-horizon games, as well as the equivalent affine variational inequality (AVI).

\begin{lemma}[{Open-Loop Nash Equilibrium via Finite-Horizon AVI, \cite[Th. 1]{benenati2025linear}}]
Assume that the augmented Riccati equations admit stabilizing solutions $(P_i,K_i)$ and $\hat P_i$, respectively, and that $\mathcal U_T(x_0) = \{ D \mathbf u + d \le 0 \}$ is nonempty. Let $\mathbb X_f \subseteq \mathcal X $ be forward-invariant for $x[t+1] = A_{\mathrm{cl}} x[t]$ (i.e., once the state enters $\mathbb X_f$, the stabilizing feedback keeps it feasible for all future times). Then any $\mathbf u^\star \in \mathcal U_T(x_0)$ solving the affine variational inequality
\begin{align}\label{eq: AVI problem}
    \mc P(x_0): \quad \langle M \mathbf u^\star + q, \mathbf u - \mathbf u^\star \rangle \ge 0, \,\, \forall \mathbf u \in \mathcal U_T(x_0), \quad 
\end{align}
is a finite-horizon Nash equilibrium for the cost $J_i^T$. Moreover, letting $x_T = \zeta(T,x_0,\mathbf u^\star) \in \mathbb X_f$, the infinite-horizon sequence
\begin{equation}
    u_i[t] = 
    \begin{cases}
    u_i^\star[t], & t<T,\\
    K_i A_{\mathrm{cl}}^{t-T} x_T, & t \ge T
    \end{cases}
\end{equation}
is an open-loop Nash equilibrium for the cost function $J_i^\infty$.
\end{lemma}


\section{Methodology and Convergence Analysis}\label{sec: algorithm and convergence}


Let us consider that the feasible set is polyhedral,
i.e., $\mathcal C := \{u\in\mathbb R^n \mid Du+d\le0\}$, and the operator $F(u)=Mu+q$ is affine and strongly monotone, hence
the variational inequality $\mathrm{AVI}(\mathcal C,M,q)$
admits a unique solution \cite[Thm.~2.3.3]{facchinei2003finite}.
As we see later, strong monotonicity is also be instrumental in establishing the well-posedness of the Newton step.

A classical characterization of the solution of $\mathrm{AVI}(\mathcal C,M,q)$ is provided by its Karush--Kuhn--Tucker (KKT) conditions. In our affine case, solving the variational inequality is equivalent to finding $(x,\lambda)\in\mathbb R^n\times\mathbb R^m$ such that \cite[Proposition~1.2.1]{facchinei2003finite}
\begin{subequations}\label{eq:KKT-newton}
\begin{align}
Mu+q-D^\top\lambda &= 0 \label{eq:KKT-newton-a}\\
Du+d &\le 0 \label{eq:KKT-newton-b}\\
\lambda &\ge 0 \label{eq:KKT-newton-c}\\
\lambda^\top(Du+d) &= 0. \label{eq:KKT-newton-d}
\end{align}
\end{subequations}

The complementarity conditions
\eqref{eq:KKT-newton-b}--\eqref{eq:KKT-newton-d} are inherently
nonsmooth, which prevents the direct application of Newton-type methods. 
One can however reformulate these conditions with a 
\emph{nonlinear complementarity problem (NCP)}, which transforms complementarity constraints into a system of smooth equations while preserving exact equivalence. Among the most widely
used NCP-functions is the Fischer--Burmeister function
\cite{andreas1995local} 
\begin{align*}
\phi(a,b) := \sqrt{a^2+b^2} - a - b,
\end{align*}
which satisfies the fundamental equivalence
\begin{align*}
\phi(a,b)=0 \quad \Longleftrightarrow \quad a \ge 0, \; b \ge 0, \; ab = 0,
\end{align*}
and allows one to equivalently reformulate the nonsmooth complementarity conditions as $\phi(a,b)=0$. To improve numerical stability and ensure differentiability everywhere, we introduce a smoothing parameter $\mu>0$ and consider the smoothed function
\begin{align*}
\phi_\mu(a,b) := \sqrt{a^2+b^2+\mu^2} - a - b,
\end{align*}
which is continuously differentiable for all $(a,b)\in\mathbb R^2$ and
satisfies $\phi_\mu(a,b)\to\phi(a,b)$ as $\mu\to 0$ 
\cite{chen1997smooth,kanzow1996some}.

Applying this transformation to the complementarity constraints in
\eqref{eq:KKT-newton} leads to the smooth system
\begin{align}\label{eq: NCP-function}
\Phi_\mu(u,\lambda) :=
\begin{bmatrix}
Mu+q-D^\top \lambda\\[1mm]
\phi_\mu(Du+d,\lambda)
\end{bmatrix} = 0,
\end{align}
where $\phi_\mu(Du+d,\lambda)$ is understood componentwise.  
For each component $i=1,\dots,m$, we define the partial derivatives of
$\phi_\mu$ with respect to its first and second arguments as
\begin{align*}
g_\mu^i(u,\lambda) &:= \frac{\partial \phi_\mu(a_i,b_i)}{\partial a_i} 
= \frac{a_i}{\sqrt{a_i^2+b_i^2+\mu^2}} - 1,\\
h_\mu^i(u,\lambda) &:= \frac{\partial \phi_\mu(a_i,b_i)}{\partial b_i} 
= \frac{b_i}{\sqrt{a_i^2+b_i^2+\mu^2}} - 1,
\end{align*}
with $a_i=(Du+d)_i$ and $b_i = \lambda_i$. The diagonal matrices
\begin{align*}
G_\mu := \mathrm{diag}(g_\mu^1,\dots,g_\mu^m),\,\, H_\mu := \mathrm{diag}(h_\mu^1,\dots,h_\mu^m),
\end{align*}
collect these derivatives for all components. With this notation, the Jacobian of
$\Phi_\mu$ admits the block structure
\begin{align}\label{eq: diff-NCP-function}
\nabla \Phi_\mu(u,\lambda) =
\begin{bmatrix}
M & -D^\top\\
G_\mu D & H_\mu
\end{bmatrix}.
\end{align}

Given an iterate $(u^k,\lambda^k)$, the Newton direction
$(\Delta u^k, \Delta \lambda^k)$ is computed by solving the linear system
\begin{equation}\label{eq:newton-step}
\nabla \Phi_\mu(u^k,\lambda^k)
\begin{bmatrix}
\Delta u^k\\
\Delta \lambda^k
\end{bmatrix} 
= -\Phi_\mu(u^k,\lambda^k).
\end{equation}

The nonsingularity of $\nabla\Phi_\mu$ is guaranteed for $\mu>0$ due to
the smoothing of $G_\mu$ and $H_\mu$, and the strong monotonicity of \(M\). Specifically, $H_\mu$ has
strictly negative diagonal entries and is therefore invertible. Eliminating
$\Delta\lambda^k$ from \eqref{eq:newton-step} yields the reduced system
$\bigl(M + D^\top H_\mu^{-1} G_\mu D\bigr) \Delta u^k = \tau^k$
for some right-hand side $\tau^k$. Since $G_\mu$ and $H_\mu$ are diagonal with
negative entries, $D^\top H_\mu^{-1} G_\mu D$ is positive semidefinite, and
strong monotonicity of $M$ implies that the reduced system, and therefore the
full Jacobian $\nabla \Phi_\mu$, is nonsingular for all $(u,\lambda)$; see
\cite[Ch.~7.3]{facchinei2003finite}, \cite{chen1997smooth,kanzow1996some}.


We are now ready to show fast local convergence and global convergence of the proposed Newton method (Algorithm~\ref{alg: Newton-step}).

\begin{theorem}[Local superlinear convergence of smoothed Newton]\label{Th: local convergence}
Let $(u_\mu^\star,\lambda_\mu^\star)$ be the unique solution of the
smoothed system $\Phi_\mu(u,\lambda)=0$ for $\mu>0$, and assume that the
Jacobian $\nabla \Phi_\mu(u_\mu^\star,\lambda_\mu^\star)$ is nonsingular.
Then, for iterates $(u^k,\lambda^k)$ sufficiently close to $(u_\mu^\star,\lambda_\mu^\star)$,
the Newton method defined by \eqref{eq:newton-step} converges superlinarly:
\begin{align*}
\| (u^{k+1},\lambda^{k+1}) - (u_\mu^\star,\lambda_\mu^\star) \|
\le C \| (u^{k},\lambda^{k}) - (u_\mu^\star,\lambda_\mu^\star) \|^2,
\end{align*}
for some constant $0< C < 1$ independent of $k$; see, e.g., \cite[Thm.~7.2.5]{facchinei2003finite} and \cite[Thm. 6C.1 \& 6E.2]{dontchev2009implicit}.
\end{theorem}


\begin{theorem}[Global convergence of line-search Newton]\label{thm: global conv}
Let $\Phi_\mu:\mathbb{R}^{n}\times\mathbb{R}^{m}\to\mathbb{R}^{n+m}$ be continuously
differentiable with nonsingular Jacobian $\nabla \Phi_\mu(u_\mu^\star,\lambda_\mu^\star)$ at
its unique solution $(u_\mu^\star,\lambda_\mu^\star)$. Let $(u^k,\lambda^k)$ be the
iterates generated by the Newton direction
\begin{align*}
\nabla \Phi_\mu(u^k,\lambda^k)
\begin{bmatrix}
\Delta u^k\\
\Delta \lambda^k
\end{bmatrix} = -\alpha_k \Phi_\mu(u^k,\lambda^k),
\end{align*}
with stepsize $\alpha_k \in (0,1]$ determined by a merit-function linesearch $\Psi(u,\lambda)$ satisfying a standard Armijo condition \cite{Armijo1966}. Then:
\begin{enumerate}
\item[(i)] From any initial point $(u^0,\lambda^0)$, the iterates converge to
$(u_\mu^\star,\lambda_\mu^\star)$.
\item[(ii)] Once the iterates enter a neighborhood of $(u_\mu^\star,\lambda_\mu^\star)$,
the linesearch accepts unit stepsize ($\alpha_k = 1$) and the iterates
converge superlinearly to the solution.
\end{enumerate}

We refer to \cite[Ch.~6]{dennis1996numerical} and \cite[Ch.~8.3]{facchinei2003finite} for technical proofs.

\end{theorem}


\begin{algorithm}[!h]
\caption{Newton Step for Smoothed System $\Phi_\mu(u,\lambda)=0$}
\label{alg: Newton-step}
\begin{algorithmic}[1]
\REQUIRE Initial point $(u^0,\lambda^0)$, smoothing parameter $\mu>0$
\STATE Evaluate $\Phi_\mu(u^k,\lambda^k)$ \eqref{eq: NCP-function}, $\nabla \Phi_\mu(u^k,\lambda^k)$~\eqref{eq: diff-NCP-function}, and the Newton direction $d_k$
\[
d_k = (\Delta u^k,\Delta \lambda^k) : = - \Big(\nabla \Phi_\mu(u^k,\lambda^k)\Big)^{-1}\Phi_\mu(u^k,\lambda^k),
\]
\STATE Compute the stepsize $\alpha_k$~(Algorithm~\ref{alg: Armijo})
\STATE Update the next iteration:
\[
\begin{bmatrix} u^{k+1}\\ \lambda^{k+1}\end{bmatrix} = \begin{bmatrix} u^{k}\\ \lambda^{k}\end{bmatrix} + \alpha_k d_k
\]
\end{algorithmic}
\end{algorithm}
\begin{algorithm}[!h]
\caption{Armijo linesearch for Newton Direction}
\label{alg: Armijo}
\begin{algorithmic}[1]
\REQUIRE Newton direction $d_k = (\Delta u^k,\Delta \lambda^k)$, current iterate $(u^k,\lambda^k)$,
merit function $\Psi(u,\lambda) = \frac12 \|\Phi_\mu(u,\lambda)\|^2$, 
parameters $c\in(0,1)$, and $\beta\in(0,1)$
\STATE Initialize step size $\alpha \gets 1$
\WHILE{$\Psi(u^k+\alpha \Delta u^k, \lambda^k+\alpha \Delta \lambda^k) 
       > \Psi(u^k,\lambda^k) + c \alpha \nabla \Psi^\top (\Delta u^k,\Delta \lambda^k)$}
    \STATE $\alpha \gets \beta \alpha$
\ENDWHILE
\STATE Return step size $\alpha$
\end{algorithmic}
\end{algorithm}

Let us emphasize that the main limitation of Newton-type methods is the computational complexity for solving the system \eqref{eq:newton-step} to determine the descent direction. In this work, we exploit the special structure of our Jacobian \eqref{eq: diff-NCP-function} to reduce the complexity from $\mathcal{O}\big((n+m)^3\big)$ to $\mathcal{O}\big( n^3 + m n^2 \big)$.

\section{Applications and Numerical Experiments}\label{sec: simulation}

We demonstrate the performance of Algorithm~\ref{alg: Newton-step} on \textit{(A) Vehicle platooning}, adapted from~\cite{shi2017distributed,benenati2025linear}, and \textit{(B) Unsignalized intersection crossing} ~\cite{baghbadorani2025douglas}. To evaluate performance, we compare our proposed algorithm (NW), including a high-efficiency implementation (Fast-NW), versus two extensively used methods: (i) Forward-Backward method (FB)~\cite{nemirovskij1983problem} and (ii) the Douglas-Rachford method (DR)~\cite{ferris1996operator,baghbadorani2025douglas}, the state-of-the-art approach for AVI problems. In all simulations, the smoothing parameter is $\mu = 10^{-6}$ and the initial points are chosen randomly. Additionally, we set the linesearch parameter to $c = 10^{-4}$ and the backtracking factor to $\beta = 0.5$ within the Armijo linesearch (Algorithm~\ref{alg: Armijo}). We terminate all iterative algorithms if $\|r_k\| \leq 10^{-4}$, where $r_k$ denotes the residual at iteration $k$, or if the number of iterations reaches $10^5$.
\subsection{Vehicle Platooning}

Let us adopt the vehicle platooning scenario as described in \cite{shi2017distributed,benenati2025linear}. We consider $N$ vehicles, where the leading vehicle (indexed by $1$) tracks a reference velocity $v^{\mathrm{ref}}$, while the subsequent agents $i \in \{2, \dots, N\}$ aim to synchronize their speed with the preceding vehicle. To maintain safety, each agent must keep an inter-vehicle distance defined by a fixed gap $d_i$ and a velocity-dependent headway $h_i v_i$. Following \cite{benenati2025linear}, the local state for agent $i > 1$ is defined as $x_i = [p_{i-1} - p_i - d_i - h_i v_i, v_{i-1} - v_i]^\top$, where $p_i$ and $v_i$ represent the position and velocity of the $i$-th agent. For the leader, the state is relative to the reference: $x_1 = [0, v^{\mathrm{ref}} - v_1]^\top$.

The agents are modeled using sampled double-integrator dynamics with a sampling period $\tau_{\text{s}} = 0.1\,\text{s}$. The discrete-time collective dynamics are represented as $x^+ = Ax + \sum B_i u_i$, where the system matrices are defined as follows:
\begin{align} \label{eq:platooning_dyn}
    \begin{split}
        A &= \text{diag}\left(\begin{bmatrix} 0 & 0 \\ 0 & 1 \end{bmatrix}, I_{N-1} \otimes \begin{bmatrix} 1 & \tau_{\text{s}} \\ 0 & 1 \end{bmatrix} \right); \\
        B_1 &=  \delta_{2}^N \otimes \begin{bmatrix} \tau_{\text{s}}^2/2 \\ \tau_{\text{s}} \end{bmatrix} - \delta_1^N \otimes \begin{bmatrix} 0 \\ \tau_{\text{s}} \end{bmatrix};\\
        B_N &= -\delta_N^N \otimes \begin{bmatrix} h_i\tau_{\text{s}} + \tau_{\text{s}}^2/2 \\ \tau_{\text{s}} \end{bmatrix};\\
        B_i &=  \delta_{i+1}^N \otimes \begin{bmatrix} \tau_{\text{s}}^2/2 \\ \tau_{\text{s}} \end{bmatrix} - \delta_i^N \otimes \begin{bmatrix} h_i\tau_{\text{s}} + \tau_{\text{s}}^2/2 \\ \tau_{\text{s}} \end{bmatrix},  \\
        &\forall i\in\{2,...,N-1\}.
    \end{split}
\end{align}

Within this framework, we enforce several critical constraints: a safety gap $p_{i-1} \geq d^{\text{min}}_i + p_i$, velocity limits $v_i \in [v_i^{\text{min}}, v_i^{\text{max}}]$, and input saturation $u_i \in [u_i^{\text{min}}, u_i^{\text{max}}]$. Because the open-loop dynamics in \eqref{eq:platooning_dyn} do not satisfy standard stabilizability assumptions, we apply the pre-stabilizing local controller $K_i^{\text{stab}} = (\delta_i^N)^{\top} \otimes [-1, -1]$ proposed in the reference study.
\begin{figure}[!h]
    \centering
    \includegraphics[width=1\linewidth]{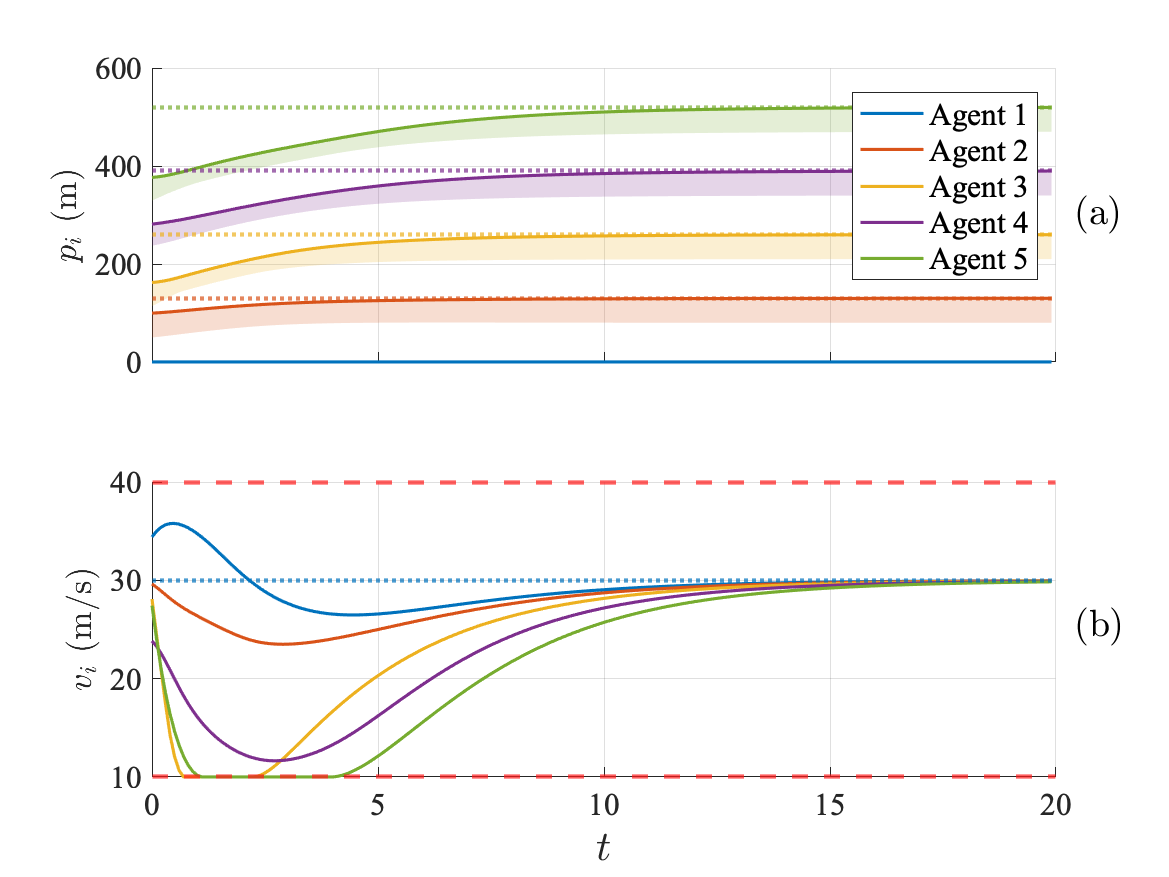}
    \caption{(a) Position of vehicles with respect to the leading agent (b) Velocity of each agent (dotted lines: reference values; red dashed lines: constraints).}
    \label{fig: platooning-traj}
\end{figure}
\begin{figure}[!h]
    \centering
    \includegraphics[width=.85\linewidth]{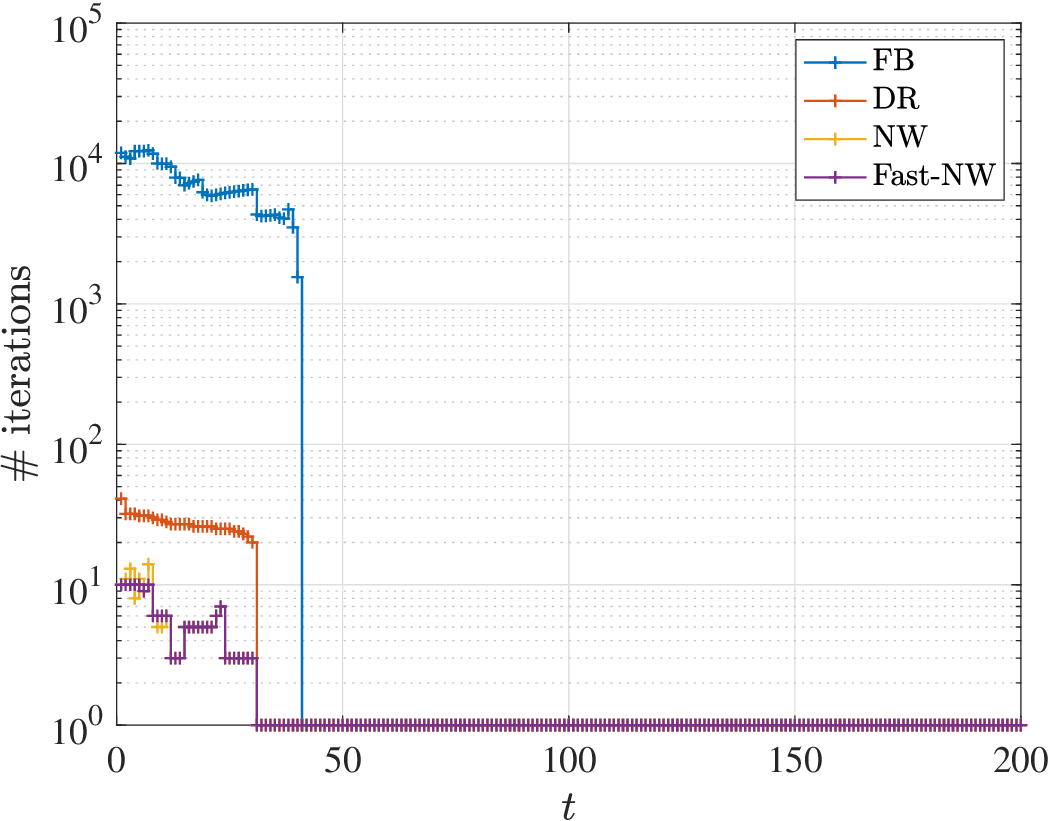}
    \caption{{Iterations required for convergence of the VI solvers.}} 
    \label{fig: platooning_iter_convergence}
\end{figure}
\begin{figure}[!h]
    \centering
    \includegraphics[width=.85\linewidth]{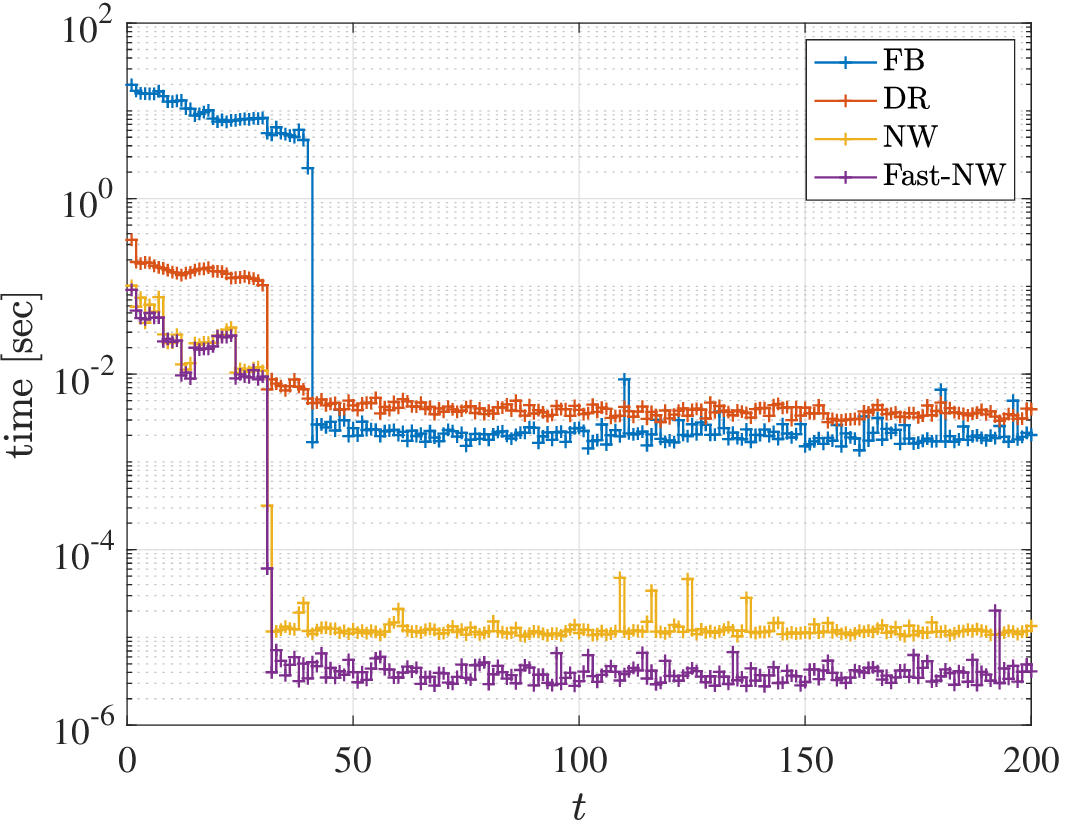}
    \caption{{Computational time for convergence of the VI solver.}} 
    \label{fig: platooning_time_convergence}
\end{figure}
In our comparative analysis, we utilize the OL-NE receding-horizon control framework with a prediction horizon $T=10$ and weights $Q_i = I, R_i = 1$. Figures~\ref{fig: platooning_iter_convergence} and \ref{fig: platooning_time_convergence} illustrate the number of iterations and the computational time required at each time step to solve the variational inequality problem. The (Fast)-Newton methods outperform all other approaches in terms of both iteration count and total execution time. Furthermore, Figure~\ref{fig: platooning-traj} shows the resulting vehicle trajectories, demonstrating that the agents achieve the desired equilibrium state while satisfying all system constraints.

\subsection{Unsignalized intersection crossing}

We consider a second scenario involving the coordination of $N=15$ autonomous vehicles crosssing an insignalized intersection, with first-come, first-served priority sequence based on the vehicles' arrival times \cite{baghbadorani2025douglas}.
In this scenario, the traffic flow is characterized by a diverse set of maneuvers, including straight crossings (e.g., North-to-South (NS), East-to-West (EW)) and various turning trajectories (e.g., Northwest (NW) or West-to-South (WS)). The primary coordination challenge lies in managing the coupling between vehicles with intersecting trajectories; we define $\chi(i)$ as the index of the preceding vehicle occupying a conflicting path for agent $i$. For agents in the set of leading vehicles $\mathcal{L}$, the objective is to track a reference speed $v^{\text{ref}}$, while followers strive to maintain a safety margin $d_i$ relative to their assigned predecessor $\chi(i)$.

To capture these dependencies, the local state $x_i$ for each agent is defined as:
\begin{equation}
    x_i = \begin{cases} 
        v^{\text{ref}} - v_i & \text{if } i \in \mathcal{L} \\
        [p_{\chi(i)} - p_i - d_i, \, v_{\chi(i)} - v_i]^\top & \text{if } i \notin \mathcal{L},
    \end{cases}
\end{equation}
where $p_i$ and $v_i$ represent longitudinal progress and velocity. The collective system is modeled as a set of double-integrators discretized with a sampling rate $\tau_{\text{s}}=0.1$s, leading to the global representation $x^+ = Ax + \sum B_i u_i$. The system matrices $A = \text{blkdiag}(A_i)_{i\in\mathcal{I}}$ and $B_i = \text{col}(B_{ij})_{j\in\mathcal{I}}$ are structured as follows to reflect the hierarchical coupling between agents:
\begin{align} \label{eq:intersection_dynamics}
    \begin{split}
        A_i &= \begin{cases} 1 & \text{if } i \in \mathcal{L} \\ \left[\begin{smallmatrix} 1 & \tau_{\text{s}} \\ 0 & 1 \end{smallmatrix}\right] & \text{if } i \notin \mathcal{L} \end{cases}, \\
        B_{ij} &= \begin{cases} 
            [\frac{\tau_{\text{s}}^2}{2}, \tau_{\text{s}}]^\top & \text{if } i = \chi(j) \\
            -[\frac{\tau_{\text{s}}^2}{2}, \tau_{\text{s}}]^\top & \text{if } i=j, i \notin \mathcal{L} \\
            -\tau_{\text{s}} & \text{if } i=j, i \in \mathcal{L} \\
            0 & \text{otherwise.}
        \end{cases}
    \end{split}
\end{align}

Our formulation ensures that an input $u_i$ from a predecessor directly influences the state error of its subsequent follower $j$. To satisfy the underlying stabilizability assumptions, we apply a decentralized pre-stabilizing law $\bar{K}_i x = -0.1 \cdot \mathds{1}^\top x_i$. Furthermore, the vehicles are subject to physical and safety constraints: minimum longitudinal gap $p_{\chi(i)} - p_i \geq d_{\text{min}}$, velocity $v_i \in [v^{\text{min}}, v^{\text{max}}]$ and input bounds $u_i \in [u^{\text{min}}, u^{\text{max}}]$.
\begin{figure}[!h]
    \centering
    \includegraphics[width=1\linewidth]{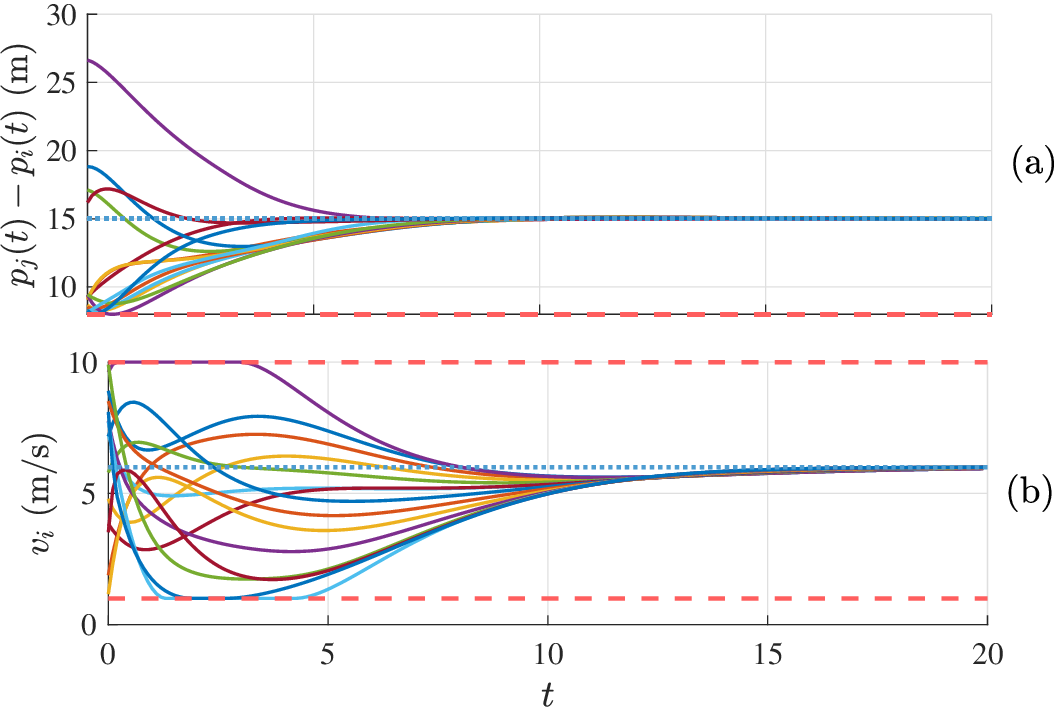}
    \caption{(a) Distance between $\chi(i)$ and $i$. (b) Velocities.}
    \label{fig: crossroad-traj}
\end{figure}

\begin{figure}[!h]
    \centering
    \includegraphics[width=.85\linewidth]{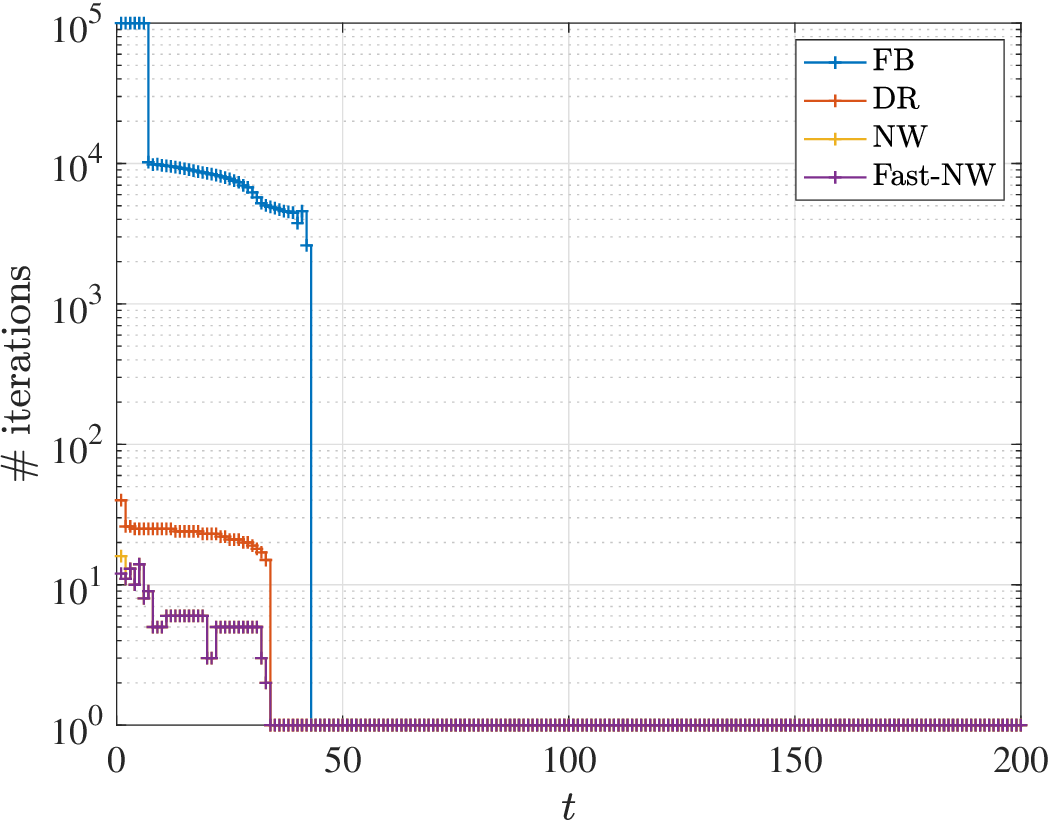}
    \caption{{Iterations required for convergence of the VI solvers.}} 
    \label{fig: crossroad_iter_convergence}
\end{figure}

\begin{figure}[!h]
    \centering
    \includegraphics[width=.85\linewidth]{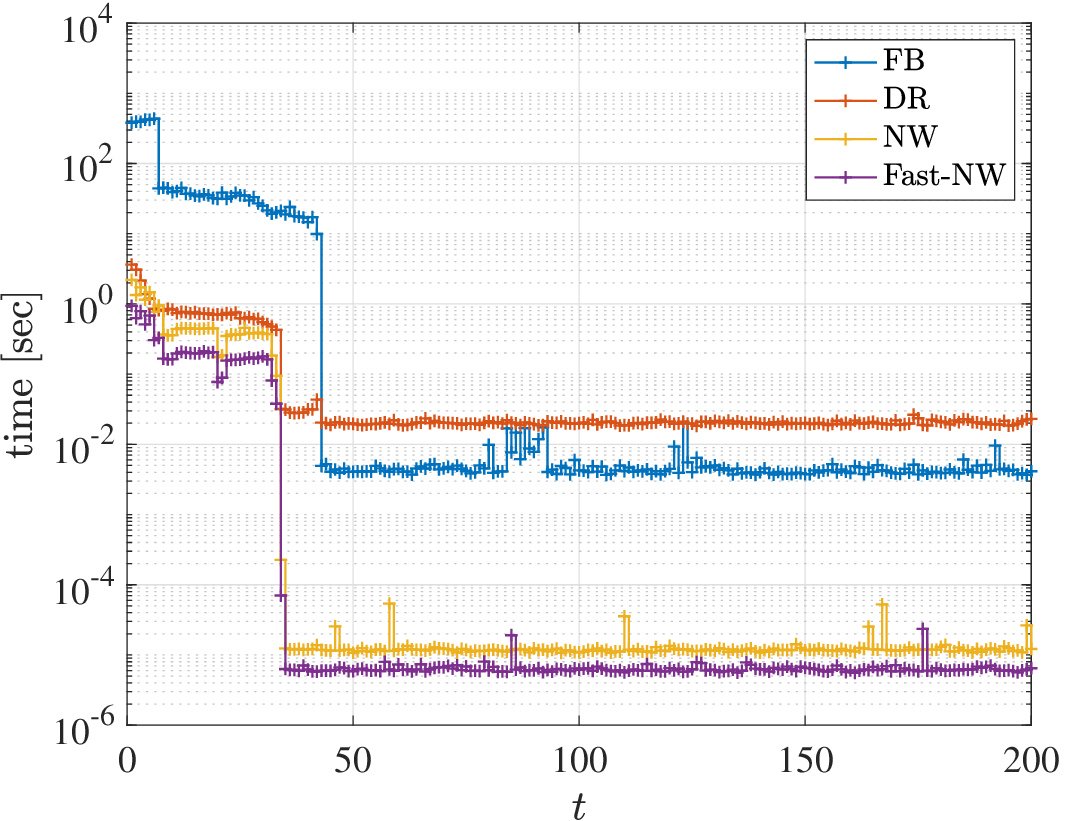}
    \caption{{Computational time for convergence of the VI solvers.}} 
    \label{fig: crossroad_time_convergence}
\end{figure}

Let us adopt weighting matrices $Q_i = I$ and $R_i = 1$, supplemented by a warm-start strategy where the input sequence is shifted at each time step. As illustrated in Figures~\ref{fig: crossroad_iter_convergence} and \ref{fig: crossroad_time_convergence}, the (Fast)-Newton approaches have superior computational efficiency, significantly outperforming standard Forward-Backward (FB) and Douglas-Rachford (DR) splitting methods in both convergence rate and execution time. Furthermore, the state trajectories in Figure~\ref{fig: crossroad-traj} confirm that the Newton-based solvers enable all vehicles to reach the reference velocity and maintain desired inter-vehicle spacing while satisfying all safety constraints.
The resulting closed-loop dynamics are illustrated in this animation video: \texttt{http://bit.ly/4aLhsSg}.


\begin{figure}[!h]
    \centering
    \includegraphics[width=1\linewidth]{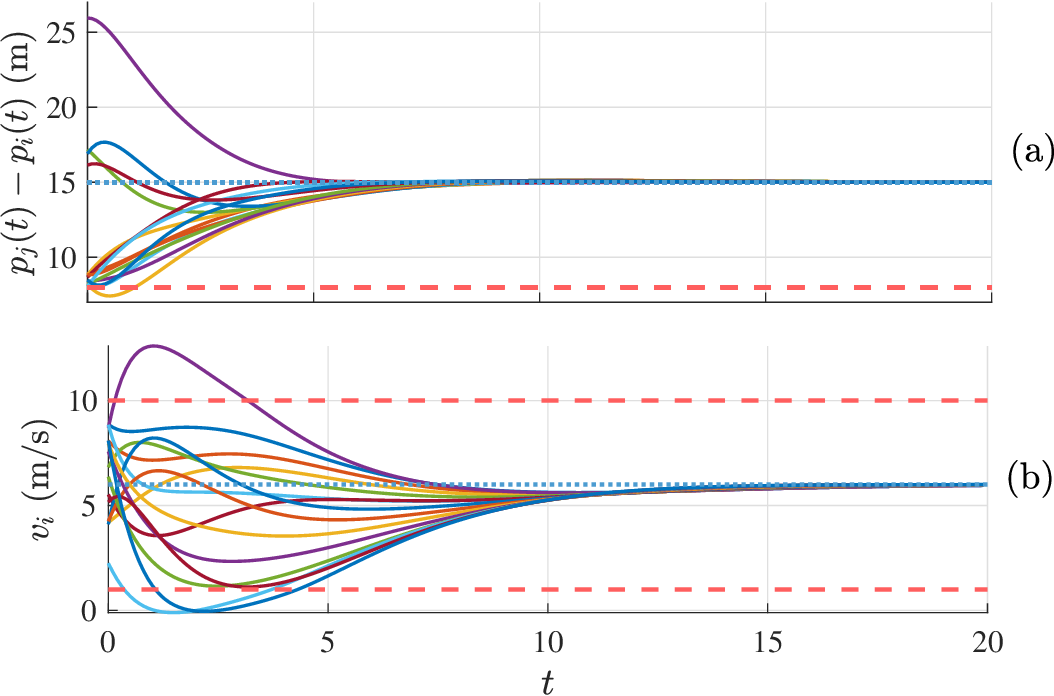}
    \caption{Forward-Backward trajectories for 10 fixed iterations. (a) Yellow trajectory violates distance constraints. (b) Violet, blue, and cyan agents do not satisfy velocity constraints.}
    \label{fig: FB_violation_iterations}
\end{figure}

\begin{figure}[!h]
    \centering
    \includegraphics[width=.9\linewidth]{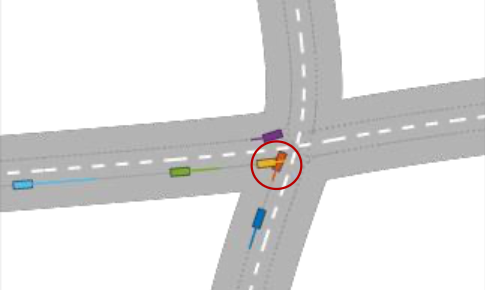}
    \caption{Douglas–Rachford trajectories for 10 fixed iterations. The orange and yellow agents collide as illustrated in this animation video: \texttt{https://bit.ly/3Oro5Sq}.}
    \label{fig: DR_collision}
\end{figure}

To further evaluate the robustness of the Newton method, we consider a restricted computational budget where the solver is limited to $10$ iterations per time step, a scenario representative of high-frequency real-time control. The resulting trajectories, plotted in Figure~\ref{fig: FB_violation_iterations}, reveal that the FB method fails to ensure feasibility within this limited budget, leading to constraint violations. In addition, Figure~\ref{fig: DR_collision} illustrates a safety distance violation and subsequent collision under the DR method, whereas our proposed Newton approach maintains a safe operating regime even under strict computational limits.


\section{Conclusion and outlook}
The Newton method allows us to solve affine variational inequalities extremely fast, thus making receding-horizon optimal control of constrained linear-quadratic dynamic games applicable to intelligent and connected autonomous vehicles.

In the future, we will investigate quasi-Newton methods, which utilize first-order information to approximate the local curvature. Furthermore, we plan to validate our receding-horizon control approach in realistic laboratory experiments.


\bibliographystyle{IEEEtran} 
\bibliography{references.bib}
\end{document}